\documentclass[12pt]{article}
\usepackage[margin=1in]{geometry}
\usepackage{url}

\usepackage[T1]{fontenc}
\usepackage{fouriernc}

\usepackage{amsmath}
\usepackage{amssymb}
\usepackage{amsthm}
\usepackage{mathrsfs}

\theoremstyle{plain}
\newtheorem{theorem}{Theorem}

\theoremstyle{definition}

\newtheorem{example}[theorem]{Example}

\theoremstyle{remark}

\title{Common factors in automatic and Sturmian sequences}
\author{
Narad Rampersad\\
Department of Mathematics and Statistics\\
University of Winnipeg\\
\url{n.rampersad@uwinnipeg.ca}\\
\and
Jeffrey Shallit\\
School of Computer Science\\
University of Waterloo\\
\url{shallit@uwaterloo.ca}}
\date{\today}

\begin{document}
\maketitle

\begin{abstract}
In this short note we show that a $k$-automatic sequence and a
Sturmian sequence cannot have arbitrarily large factors in common.
\end{abstract}

\section{Introduction}
Sturmian sequences are those given by the first differences
of sequences of the form 
$$(\lfloor n \alpha  + \beta\rfloor)_{n \geq 1},$$ where
$0 \leq \alpha, \beta < 1$ and $\alpha$ is irrational
\cite{BS}.
It is well-known that a Sturmian sequence cannot be $k$-automatic; that is,
it cannot be generated
by a finite automaton reading $n$ expressed in an integer
base $k \geq 2$. This follows from the fact that the limiting frequency
of $1$'s in a Sturmian sequence is $\alpha$, whereas if a letter in
a $k$-automatic sequence has a limiting frequency, that frequency must be
rational \cite[Thm.~6, p.~180]{C}. 

Recall that by {\it factor} of a word or sequence $\bf s$, we mean a contiguous
block of symbols $x$ inside $\bf s$.
Then a natural
question is, can a Sturmian sequence and a $k$-automatic sequence
have arbitrarily large finite factors in common?  This question is related to a
problem recently studied by Byszewski and Konieczny \cite{BK17}: they wish to
determine which \emph{generalized polynomial functions} (these are
sequences defined by expressions involving algebraic operations along
with the floor function) can be $k$-automatic\footnote{We obtained these
results in July 2016.  The result was also mentioned at the \textit{Bridges between Automatic
    Sequences, Algebra, and Number Theory} School held at the CRM in
  Montreal in April 2017, where Jakub Byszewski pointed out the
  connections to his work.  It is for this reason that we are posting
  the proof of this result.}.

We also mention the work of Tapsoba \cite{Tap95}.  Recall that the complexity of a word $\bf s$ is the function counting the number of distinct 
factors of length $n$ in $s$.   It is also well-known
that Sturmian words have the minimum possible 
complexity $n+1$ achievable by an aperiodic
infinite word.   Tapsoba shows another distinction between automatic sequences
and Sturmian words by giving a formula for the minimal complexity function of the
fixed point of an injective $k$-uniform binary morphism and comparing this to
the complexity function of Sturmian words.

Our main result is the following:

\begin{theorem}\label{main}
Let $\bf x$ be a $k$-automatic sequence and let $\bf a$ be a Sturmian
sequence.  There exists a constant $C$ (depending on $\bf x$ and $\bf a$) such
that if $\bf x$ and $\bf a$ have a factor in common of length $n$, then $n \leq
C$.
\end{theorem}

Note that this result would follow fairly easily from the frequency
results mentioned previously, \emph{if} $\bf x$ is
\emph{uniformly recurrent} (meaning that for every factor $z$ of $\bf x$
occurs infinitely often, and with bounded gap size between two consecutive occurrences).  However, unlike Sturmian sequences, automatic sequences
need not be uniformly recurrent:  consider, for example, the $2$-automatic
sequence that is the characteristic sequence of the powers of $2$.
Our proof is therefore based
on the finiteness of the \emph{$k$-kernel} of $\bf x$, along with the
\emph{uniform distribution} property of Sturmian sequences (similar
arguments have previously been used by the second author \cite{Sha96}).

\section{Proof of Theorem~\ref{main}}
\begin{proof}
Let ${\bf x}=x_0x_1 \cdots$ and ${\bf a}=a_0a_1\cdots$.  Since the factors of a Sturmian word do not depend on $\beta$, without loss of
generality, we may suppose that $\beta = 0$ (or, in other words, that $\bf a$ is a characteristic word).  Then there
exists an irrational number $\alpha$ such that $\bf a$ is defined by the
following rule:
\[
a_n = \begin{cases}1, &\text{ if } \{(n+1)\alpha\} < \alpha;\\
0, &\text{ otherwise.}\end{cases}
\]
Here $\{\cdot\}$ denotes the fractional part of a real number.

Suppose that for some $L$, the words $\bf x$ and $\bf a$ have a factor of
length $L$ in common: i.e., for some $i \leq j$
\[ x_i \cdots x_{i+L-1} = a_j \cdots a_{j+L-1}. \]
(We may assume that $i \leq j$ since $\bf a$ is recurrent, but this is not
important for what follows.)  Suppose that the $k$-kernel of $\bf x$,
\[
\{ (x_{nk^r+s})_{n \geq 0} : r \geq 0 \text{ and } 0 \leq s < k^r \},
\]
has $Q$ distinct elements.  Let $r$ satisfy $k^r > Q$.  There there
exist integers $s_1, s_2$ with $0 \leq s_1 < s_2 < k^r$ such that
\[
(x_{nk^r+s_1})_{n \geq 0} = (x_{nk^r+s_2})_{n \geq 0}.
\]

Define
\begin{align*}
d_1 &:= s_1 + j-i + 1\\
d_2 &:= s_2 + j-i + 1.
\end{align*}
For all $n$ satisfying $i \leq nk^r+s_1$ and $nk^r+s_2 \leq i+L-1$ we
have $x_{nk^r+s_1} = a_{nk^r+d_1-1}$ and $x_{nk^r+s_2} =
a_{nk^r+d_2-1}$.  Since $x_{nk^r+s_1} = x_{nk^r+s_2}$, we have $a_{nk^r+d_1-1} =
a_{nk^r+d_2-1}$.  This means that either the inequalities
\begin{equation}\label{both0}
\{(nk^r+d_1)\alpha\} < \alpha \text{ and } \{(nk^r+d_2)\alpha\} < \alpha
\end{equation}
both hold, or the inequalities
\begin{equation}\label{both1}
\{(nk^r+d_1)\alpha\} \geq \alpha \text{ and } \{(nk^r+d_2)\alpha\} \geq \alpha
\end{equation}
both hold.

If $L$ is arbitrarily large, then there exist arbitrarily large sets
$I$ of consecutive positive integers such that every $n\in I$
satisfies either \eqref{both0} or \eqref{both1}.  Without loss of
generality, suppose that $\{d_2\alpha\} > \{d_1\alpha\}$. Choose
$\epsilon>0$ such that $\epsilon < \{d_2\alpha\} - \{d_1\alpha\}$.
Note that $d_2-d_1 = s_2 - s_1$, so $\epsilon$ does not depend on $L$
(or $I$).  Since $k^r\alpha$ is irrational, if $I$ is sufficiently
large then by Kronecker's theorem (which asserts that the set of
points $\{n\alpha\}$ is dense in $(0,1)$) there exists $N \in I$ such
that
\[
\{N(k^r\alpha)+d_2\alpha\} \in [\alpha,\alpha+\epsilon].
\]
By the choice of $\epsilon$, this implies that
\[
\{N(k^r\alpha)+d_2\alpha\} \geq \alpha \text{ and }
\{N(k^r\alpha)+d_1\alpha\} < \alpha,
\]
contradicting the assumption that $N$ satisfies one of \eqref{both0}
or \eqref{both1}.  The contradiction means that $L$ must be bounded by
some constant $C$, which proves the theorem.
\end{proof}

\begin{example}
Consider the Thue-Morse word ${\bf t} = 01101001\cdots$ given by the
fixed point of the morphism $0 \rightarrow 01$ and $1 \rightarrow 10$,
and the Fibonacci word ${\bf f} = 01001010\cdots$ given by the
fixed point of $0 \rightarrow 01$ and $1 \rightarrow 0$.  The latter
is Sturmian.  The set of common factors is
\begin{multline*}
\{ \epsilon, 0, 1, 00, 01, 10, 001, 010, 100, 101, 0010, 0100, 0101, 1001, 1010, \\ 
00101, 01001, 10010, 10100,
010010, 100101, 101001, 0100101, 1010010, 10100101 \},
\end{multline*}
so $C = 8$.
\end{example}

\section*{Acknowledgment}

We thank Jean-Paul Allouche for helpful discussions.

\end{document}